\documentclass[12pt]{amsart}
\usepackage{amscd,amssymb}
\usepackage[arrow,matrix]{xy}

\theoremstyle{plain}
\newtheorem{theorem}[subsection]{Theorem}
\newtheorem{lemma}[subsection]{Lemma}
\newtheorem{prop}[subsection]{Proposition}
\newtheorem{corollary}[subsection]{Corollary}

\theoremstyle{definition}
\newtheorem{remark}[subsection]{Remark}

\newtheorem{example}[subsection]{Example}

\numberwithin{equation}{section}
\newcommand{\U} {{\mathcal U}}

\newcommand{\HC}{{\mathcal H}}
\newcommand{\F}{{\mathcal F}}
\newcommand{\G}{{\mathcal G}}

\newcommand{\LL}{{\mathcal L}}

\newcommand{\T}{{\mathcal T}}

\renewcommand{\SS}{{\mathcal S}}

\newcommand{\Z}{\mathbb{Z}}
\newcommand{\Q}{\mathbb{Q}}

\newcommand{\C}{\mathbb{C}}

\newcommand{\PP}{\mathbb{P}}
\newcommand{\HH}{\mathbb{H}}

\begin{document}
\date{}

\title [Regular functions transversal at infinity]
{Regular functions transversal at infinity     }

\author[Alexandru Dimca]{Alexandru Dimca}
\address{A.Dimca:  Laboratoire J.A. Dieudonn\'e, UMR du CNRS 6621,
                 Universit\'e de Nice-Sophia-Antipolis,
                 Parc Valrose,
                 06108 Nice Cedex 02,
                 FRANCE.}
\email
{dimca@math.unice.fr}

\author[Anatoly Libgober]{Anatoly Libgober}
\address{A.Libgober: Department of Mathematics, 
          University of Illinois at Chicago,
          851 S.Morgan, Chicago, Illinois, 60607   
                 USA.}
\email
{libgober@math.uic.edu}



\begin{abstract}
We generalize and complete some of  Maxim's
 recent results on Alexander invariants of a polynomial transversal to the hyperplane at infinity. Roughly speaking, and surprisingly, such a polynomial behaves both topologically and algebraically (e.g. in terms of the variation of MHS on the cohomology of its smooth fibers), like a homogeneous  polynomial.

\end{abstract}

\maketitle

\section{Introduction and the main results}

In the last twenty years there has been an ever increasing interest in the 
topology and geometry of polynomial functions with a 
certain good behavior at infinity, see for instance \cite{B},  \cite{GN0}, \cite{NZ}, \cite{NN}, \cite{PZ}, \cite{S}, \cite{SiT}. 
In particular the point of view of constructible sheaves was useful,
see  \cite{D}. 
An interesting problem in this area is to understand the Alexander invariants
of the complements to affine hypersurfaces defined by such polynomial 
functions. Various approaches, some algebro-geometric, 
using the superabundances 
of linear systems associated with singularities (cf. remark 
\ref{superab} in the last section), and others, more 
topological, using the monodromy representation 
were proposed (see for instance \cite{Lib}, \cite{Lib1}, \cite{KK}, \cite{DN}, \cite{O}). 
Recently, L. Maxim has considered a similar interplay 
but in a more general framework which 
includes hypersurfaces with no restrictions on singularities and 
a new, and very natural condition of good 
behavior at infinity, that we describe now.

Let $X \subset \C^{n+1}$ with $n>1$ be a reduced 
hypersurface given by an equation $f=0$.
We say that the polynomial function 
$f: \C^{n+1} \to \C$ (or the affine hypersurface $X$) is $\infty$-transversal 
if the projective closure $V$ of $X$ in $ \PP^{n+1}$ is 
transversal in the stratified sense to the hyperplane at infinity $H=
\PP^{n+1}  \setminus  \C^{n+1}$.
Consider the affine complement $M_X= \C^{n+1} \setminus X$, 
and denote by $M_X ^c$ its {\it infinite cyclic covering} corresponding 
to the kernel of the homomorphism
$$f_{\sharp}:\pi _1(M_X) \to \pi _1(\C^*)=\Z$$
induced by $f$ and sending a class of a loop into its linking number with $X$.

Then, for any positive integer $k$, the homology group 
$H_k(M_X ^c,K)$, regarded as a module 
over the principal ideal
 domain $\Lambda_K =K[t,t^{-1}]$ with $K=\Q$ or $K=\C$,
is called the $k$-th
 Alexander module of the hypersurface 
$X$, see \cite{Lib}, \cite{DN}. When this module is torsion, we denote by 
$\Delta _k(t)$ the corresponding  $k$-th  Alexander polynomial of $X$
(i.e. the $\Lambda_K$-order of $H_k(M_X ^c,K)$).

 With this notation, one of the main results in \cite{M} 
can be stated as follows.

\begin{theorem}  \label{thm1}
Assume that  $f: \C^{n+1} \to \C$ is $\infty$-transversal.
Then, for $k<n+1$, the  Alexander modules  $H_k(M_X ^c,K)$ 
  of the hypersurface $X$ are torsion  semisimple $\Lambda_K  $-modules
 which are annihilated by $t^d-1$.

\end{theorem}

Since $M_X ^c$ is an $(n+1)$-dimensional CW complex, 
one has $H_k(M_X ^c,K)=0$  for $k>n+1$, while
  $H_{n+1}(M_X ^c,K)$ is free. In this sense, the above result is optimal.
To get a flavor of the second main result in  \cite{M} describing 
the relationship between the orders of the Alexander modules 
and the singularities of $X$, see Proposition
\ref{prop3} below.

\bigskip

Now we describe the more general setting of our paper.
Let $W'=W_0' \cup ... \cup W_m'$ be a hypersurface arrangement in $ \PP^{N}$ for $N>1$. 
 Let $d_j$ denote the degree of $W_j'$ and let $g_j=0$ be a reduced defining equation for $W_j'$ in 
$ \PP^{N}$. Let $Z \subset \PP^{N}$ be a smooth complete intersection of dimension $n+1>1$ which is not contained in $W'$ and let $W_j= W_j' \cap Z$ for $j=0,...,m$ be the corresponding hypersurface in $Z$ considered as subscheme
defined by the principal ideal generated by $g_j$.
Let $W=W_0 \cup ... \cup W_m$ denote
the corresponding hypersurface arrangement in $Z$.
We assume troughout in this paper that the following hold.

\bigskip

\noindent ({\bf H1}) All the hypersurfaces $W_j$ are distinct, reduced and irreducible; moreover  $W_0$ is smooth. 

\bigskip

\noindent ({\bf H2}   ) The hypersurface   $W_0$ is transverse in the stratified sense to $V=W_1 \cup ... \cup W_m$, i.e. if $\SS$ is a Whitney regular stratification of $V$, then  $W_0$ is transverse to any stratum $S \in \SS$.

\bigskip

 The complement $U= Z \setminus W_0$ is a smooth affine variety.
We consider the hypersurface $X= U \cap V$ in $U$ and 
its complement $M_X= U \setminus X $.
Note that $M_X=M_W$, where $M_W = Z \setminus W $. 
We use both notations, each one being related to the point of view 
(affine or projective) that we wish to emphasize.

Recall that the construction of the Alexander modules and 
polynomials was generalized in the obvious way
in  \cite{DN} to the case when $ \C^{n+1}$ is replaced by 
a smooth affine variety $U$. 
The first result is new even in the special situation considered in \cite{M}.

\begin{theorem}  \label{thm2}

Assume that $d_0$ divides the sum $\sum_{j=1,m}d_j$, say $d d_0=\sum_{j=1,m}d_j$. Then 
one has the following.

 \medskip

\noindent (i) The function
$f:U \to \C$ given by
$$f(x)= \frac{g_1(x)...g_m(x)}{g_0(x)^d}$$
is a well-defined regular function on $U$ whose generic fiber $F$ is connected.

  \medskip

\noindent (ii) The restriction $f^*: M_X \to \C^*$ of $f$ outside the zero fiber
$X$ has only isolated singularities. The affine variety $U$ has the homotopy type of a space 
obtained from $X$ by adding a number of $n$-cells equal to the sum of the Milnor numbers
of the singularities of  $f^*$.
\end{theorem}

Note that we need the connectedness of $F$ since this is one of the general assumptions made in
\cite{DN}. The second claim shows that a mapping transversal at infinity behaves like an $M_0$-tame polynomial, see \cite{D2} for the definition and the properties of $M_0$-tame polynomials. These two classes of
mappings are however distinct, e.g. the defining equation of an essential affine hyperplane arrangement
is always $M_0$-tame, but the transversality at infinity may well fail for it.

\noindent The next result says roughly that an $\infty$-transversal polynomial behaves as a homogeneous polynomial up-to (co)homology of degree $n-1$. In these degrees, the determination of the  Alexander polynomial of $X$ in $U$ is reduced to the simpler problem of computing a monodromy operator.

\begin{corollary}  \label{cor1}
With the assumption in Theorem  \ref{thm2}, the following hold.

 \medskip

\noindent (i) Let $\iota:\C^* \to\C$ be the inclusion. Then, $R^0f_*\Q_U= \Q_{\C}$ and,  for each  $0<k<n$ there is a $\Q$-local system $\LL_k$ on  $\C^*$ such that 
$$R^kf_*\Q_U=\iota   _!\LL_k.$$
 In particular, for each $0<k<n$, the monodromy operators of $f$ at the origin $T_0^k$ and at infinity 
$T_{\infty}^k$ acting on
$H^k(F,\Q)$  coincide and the above local system $\LL_k$ is precisely the local system corresponding to this automorphism of $H^k(F,\Q)$.

 \medskip
   
\noindent (ii) There is a natural morphism $H^k(M_W^c,\Q)  \to  H^k(F,\Q)$ which is an  isomorphism  for $k<n$
and a  monomorphism  for $k=n$, and
 which is compatible with the obvious actions.  In particular, the associated characteristic polynomial
$$det(tId -T_0^k)=det(tId -T_{\infty}^k  )$$
coincides to the $k$-th
 Alexander polynomial $\Delta _k(X)(t)$ of $X$ in $U$  for $k<n$ and  $\Delta _n(X)(t)$ divides the G.C.D.$(det(tId -T_{\infty}^n), det(tId -T_{0}^n))$.

\end{corollary}

The next result can be regarded as similar to some results in  \cite{CDO}, \cite{Li4} and  \cite{DLi}.
Indeed, in all these results, control over the singularities of $W$ along just one of its irreducible components (in our case along $W_0$) implies that certain local systems on the complement $M_W$ are non-resonant. See  \cite{D}, p. 218 for a discussion in the case of hyperplane arrangements.

\begin{theorem}  \label{thm3}

Let $g=g_0...g_m=0$ be the equation of the  hypersurface arrangement $W$ in $Z$ and let $F(g)$ be the corresponding
global Milnor fiber given by $g=1$ in the cone $CZ$ over $Z$.
Then 
$$H^j(F(g),\Q)=H^j(M_W,\Q)$$
 for all $j<n+1$.
In other words, the action of the monodromy on $H^j(F(g),\Q)$ is trivial  for all $j<n+1$.

\end{theorem}

The main result of our paper is the following extension of Maxim's result stated in \ref{thm1} to our more general setting described above.

\begin{theorem} \label{mainthm}

Assume that $d_0$ divides the sum $\sum_{j=1,m}d_j$, say $d d_0=\sum_{j=1,m}d_j$. Then 
the following hold.

 \medskip

\noindent (i) The Alexander modules  $H_k(M_X ^c,\Q)$ 
  of the hypersurface $X$ in $U$ are torsion  semisimple $\Lambda _{\Q}$-modules
 which are annihilated by $t^d-1$  for $k<n+1$.
 
\medskip

\noindent (ii) For $k<n+1$, the Alexander module  $H^k(M_X ^c,\Q)$ 
  of the hypersurface $X$ in $U$ has a canonical mixed Hodge structure, compatible with the action of $\Lambda _{\Q}$, i.e.
$t:H^k(M_X ^c,\Q) \to H^k(M_X ^c,\Q)$ is a MHS isomorphism. Moreover, there is an epimorphism of MHS
$p_d^*:H^k(M_X ^d,\Q) \to H^k(M_X ^c,\Q)$, where  $M_X ^d$ is the $d$-cyclic covering $M_X $ and $p_d:M_X^c \to
M_X ^d$ is the induced infinite cyclic covering.
 \end{theorem}

 Dually, for $k<n+1$, the Alexander module  $H_k(M_X ^c,\Q)$ 
  of the hypersurface $X$ in $U$ has a canonical mixed Hodge structure, compatible with the natural embedding of  $H_k(M_X ^c,\Q)$ as a subspace in the homology $H_k(M_X ^d,\Q)$.

The proof of the second claim in the above theorem, given in the last section, yields also the following consequence, saying that our  regular function $f$ behaves like a homogeneous polynomial.

\begin{corollary} \label{cor2}

With the above assumptions, the MHS on the cohomology  $H^k(F_s,\Q)$ of a smooth fiber $F_s$ of $f$ is independent of $s$ for $k<n$. In this range, the isomorphism  $H^k(M_X ^c,\Q) \to H^k(F_s,\Q)$ given by Corollary  \ref{cor1} (ii) is an
isomorphism of MHS.

\end{corollary}

MHS on Alexander invariants have already been considered in the case
of hypersurfaces with isolated singularities in \cite{Lib1}
(case of plane curves considered in also in \cite{KK}). 
The above relation of this MHS to the one on the cohomology groups 
$H^k(F_s,\Q)$ is new.
Notice that Corollary 1 in \cite{KK}, combined with the main result in
 \cite{D3} and Theorem 2.10.(ii) in  \cite{DN}, yields the following.

\begin{corollary} \label{cor3}

Let $f:\C^2 \to \C$ be a polynomial function such that $X=f^{-1}(0)$ is reduced and connected and the general fiber $F$ of $f$ is connected. Then the action of $t$ on $H_1(M_X ^c,\Q)$ is semisimple.

\end{corollary}
No example seems to be known where the action of $t$  on some $H_k(M_X ^c,\Q)$ is not semisimple. On the other hand, it is easy to find examples, even for 
$f:\C^2 \to \C$, where the monodromy at infinity operator $T_1^{\infty}$ is not semisimple, see \ref{ex2}.

\bigskip

Note that, though in some important cases, see for instance \cite{Lib}, the Hurewicz theorem 
gives the identification: $H_n(M_X^c,{\Z})=
\pi_n(M_X)$, the existence of a mixed Hodge structure
on the latter cannot be deduced for example from \cite{Mor} since 
loc.cit. considers only the situation when the  
action of the fundamental group  on the homotopy groups is nilpotent 
which in general is not the case for $\pi_n(M_X)$
and of course $M_X^c$ is not quasi-projective 
in general. 

\bigskip
 
The  proofs we propose below use various techniques. 
Theorem \ref{thm4} in section 2 is the main topological 
results and is established via non-proper Morse theory as developed 
by Hamm \cite{Hamm} and Dimca-Papadima \cite{DP}. The first proof of (a special case of) the first claim in Theorem  \ref{mainthm} in section 4 is based on
a version of Lefschetz hyperplane section theorem due to Goreski-MacPherson and based on stratified Morse theory.

The proofs in section 3 are based on Theorem 4.2 in \cite{DN} (which relates Alexander modules to the cohomology of a class of rank one local systems
on the complement $M_W$) and on a general idea of getting vanishing results
 via perverse sheaves (based on Artin's vanishing Theorem) introduced in  \cite{CDO}  and developped in  \cite{D}, Chapter 6. 

Finally, the proofs in the last section use the existence of a Leray spectral sequence
of a regular mapping in the category of mixed Hodge structures (MHS for short) for which we refer to
M. Saito \cite{Sa1},  \cite{Sa2} and 
 \cite{Sa3}. To show the independence of the MHS on the Alexander module  $H^k(M_X ^c,\Q)$ on the choice of a generic fiber of $f$, we use a result by Steenbrink-Zucker on the MHS on the subspace of invariant cocycles, see
\cite{SZ}.

\section{Topology of regular functions transversal at infinity}

The following easy remark is used repeatedly in the sequel. The proof is left to the reader.

\begin{lemma}  \label{lem1}
If the hypersurface $V$ in  $Z$  has a positive
dimensional singular locus, i.e. $dim V_{sing} >0$, and $W_0$ is transversal to $V$, then 
$$dim V_{sing} =dim (V_{sing} \cap W_0 )+1.$$
 In particular, the singular locus $V_{sing} $ cannot be contained in $W_0$.
\end{lemma}

Now we start the proof of  Theorem \ref{thm2}. In order to establish the first claim,
note that the closure ${\overline F}$ of $F$ is a general member of the pencil
$$g_1(x)...g_m(x)-tg_0(x)^d=0.$$
As such, it is smooth outside the base locus given by $g_1(x)...g_m(x)=g_0(x)=0$. 

If $d=1$, then for $t$ large the above equation gives a smooth hypersurface on $Z$, hence a smooth complete intersection in 
$\PP^N$ of dimension $n>0$, hence
an irreducible variety.

For $d>2$, a closer look shows that a singular point is located either at a point where at least two of the
polynomials $g_j$ for $1 \leq j  \leq n$ vanish, or at a singular point on one of the hypersurfaces
$W_j$ for $1 \leq j  \leq n$. It follows essentially by Lemma \ref{lem1}  that $codim Sing({\overline F}) \geq 3$, hence
 ${\overline F}$ is irreducible in this case as well. This implies that $F$ is connected.

The second claim is more involved. Fix a Whitney regular stratification
$\SS$ for the pair $(Z,V)$ such that $W_0$ is transverse to $\SS$. Let $\SS '$ be the
induced  Whitney regular stratification of $CZ$, the cone over $Z$,
whose strata are either the origin, or the pull-back of strata of $\SS$
under the projection $p:CZ \setminus \{0\} \to Z $.
Then the function $h=g_1 \cdots g_m:CZ \to \C$ is stratified by
the stratifications  $\SS '$ on $CZ$ and $\T=\{\C^*,  \{0\}\}$ on $\C$, i. e.
$h$ maps submersively strata of $\SS '$ onto strata of  $\T$.
Using Theorem 4.2.1 in \cite {BMM}, it follows that the stratification
 $\SS '$ satisfies the Thom condition $(a_h)$.

Let $F_0=\{x \in CZ;~~g_0(x)=1~\}$ by the global Milnor fiber of $g_0$ regarded as a function germ on the isolated CI singularity $(CZ,0)$. Since $W_0$ is smooth, it follows that $CW_0$ is an  isolated CI singularity and hence $F_0$ has the homotopy type of a bouquet of $(n+1)$-dimensional spheres.
Let $\Gamma (h,g_0)$ be the closure of the set of points $x \in (CZ \setminus CV)$ such that the differentials $d_xh$ and $d_xg_0$ are linearly dependant. Here and in the sequence we regard $h$ and $g_0$ as regular functions on the cone $CZ$, in particular
we have $Kerd_xh \subset T_xCZ$ for any $x  \in CZ \setminus \{0\}$. Then $\Gamma (h,g_0)$ is the {\it polar curve}
of the pair of functions $(h,g_0)$. To proceed, we need the following key technical result.

\begin{theorem} \label{thm4} With the above notation, the following hold.

\medskip

\noindent (i) $dim \Gamma (h,g_0) \le 1$.

\medskip

\noindent (ii) The set $\Sigma_1$ of the singularities of the restriction of the polynomial $h$ to $F_0 \setminus CV$ is finite.

\medskip

\noindent (iii) For any $t \in S^1$, the unit circle in $\C$, consider the pencil of intersections $(Z_{s,t})_{s \in \C}$
given by
$$ Z_{s,t}=CZ \cap  \{  g_0=s \}\cap  \{  h=t \}.$$
This pencil contains finitely many singular members, and each of them has only isolated singularities. Any intersection $Z_{0,t}$ is smooth.

\medskip

\noindent (iv) $F_0$ has the homotopy type of a space  obtained from $F_0 \cap CV$ by adding $(n+1)$-cells. More precisely,
for each critical value $b \in h(\Sigma_1)$ and each small closed disc $D_b$ centered at $b$, the tube $h^{-1}(D_b)$
 has the homotopy type of a space  obtained from $ h^{-1}(c)  $ for $c \in \partial D_b$  by adding a number $(n+1)$-cells
equal to the sum of the Milnor numbers of the singularities of  $ h^{-1}(b)  $.

\end{theorem}

\proof Note first that  $\Gamma (h,g_0)$ is $\C^*$-invariant. Hence, if $dim \Gamma (h,g_0) \le 1$, then $\Gamma (h,g_0)$ may be the empty set, the origin or
a finite set of lines in $CZ$ passing through the origin.

Assume that contrary to {\it (i)} one 
has $dim \Gamma (h,g_0)>1$. Then its image in $Z$ has a positive dimension 
and hence there exist a curve $C$ on $Z$ 
along which the differentials $d_xh$ and $d_xg_0$ are linearly dependant.
Let $p$ be a point in the non-empty intersection $C \cap V$. It follows that the line $L_p$ in $\C^{N+1}$ associated to $p$
is contained in $CZ$ and that $h$ vanishes along this line. The chain rule implies that $g_0$ has a zero derivative along $L_p$, hence $g_0|L_p$is constant. Since $g_0$ is a homogeneous polynomial and the line $L_p$ passes through the origin, this constant is zero, i.e.
$g_0$ vanishes along $L_p$. Therefore  $p \in W_0 \cap V$. 
If $p$ is a smooth point on $V$, this contradicts already the transversality $W_0 \pitchfork V$.
If not, let $S \in \SS$ be the stratum containing $p$.  $W_0 \pitchfork \SS$ implies that $dimS>0$.
Let $q \in L_p$ be any nonzero vector, and let $\gamma (t)$ be an analytic curve such that $\gamma (0)=q$
and $\gamma (t) \in \Gamma (h,g_0) \setminus CV$ for $0<|t|<\epsilon.$ Hence for $t \ne 0$, $h(\gamma (t)) \ne 0$
and hence $Kerd_{\gamma (t)}h=Kerd_{\gamma (t)}g_0$. Passing to the limit for $t \to 0$ we get
$$T=limKerd_{\gamma (t)}h=limKerd_{\gamma (t)}g_0=T_q(CW_0).$$
On the other hand, the Thom condition $(a_h)$ implies
$$T \supset T_qS'=T_q(CS).$$
This implies $T_pW_0 \supset T_pS$, in contradiction to  $W_0 \pitchfork \SS$. 
The above argument shows that $dim\Gamma (h,g_0) \le 1$ and hence completes the proof of (i).

\medskip

To prove (ii), just note that $d_qh|T_qF_0=0$ for some point $q \in F_0 \setminus CV$ implies
$q \in \Gamma (h,g_0)$. Since any line through the origin intersects $F_0$ in at most $d_0$ points, the claim (ii) follows.

The last claim of (iii) is clear by homogeneity. The rest is based on the fact that 
 any line through the origin intersects $g=t$ in finitely many points.

To prove (iv) we use the same approach as in the proof of Theorem 3 in \cite {DP}, based on 
Proposition 11 in loc.cit.. Namely, we start by setting  $A=F_0$ and $f_1=h$ and construct inductively the other polynomials 
$f_2$, ...,$f_{N+1}$ to be generic homogeneous polynomials of degree $d_0$ as in loc.cit. p.485 (where generic linear forms are used for the same purpose). For more details on the non-proper Morse theory used here we refer to Hamm
\cite{Hamm}.

\endproof

We continue now the proof of the second claim in  Theorem \ref{thm2}.
There is a cyclic covering $F_0 \to U$ of order $d_0$ which restrict to a similar covering
$$p: F_0 \setminus CV \to U \setminus X $$
satisfying $f=h \circ p$. Using this and the claim (ii) above we get
that the restriction $f^*: U \setminus X \to \C^*$ of $f$ has only isolated singularities. Let $G$ be the cyclic group of order  $d_0$. Then $G$ acts on $F_0$
as the monodromy group of the function $g_0$, i.e. the group spanned by the monodromy homeomorphism
$x \mapsto \kappa \cdot x  $ with $\kappa=exp {{2 \pi  \sqrt {-1}}\over d_0}$. Since $d_0|d$, the function $h$ is $G$-invariant.
Note that the above construction of $F_0$  from $F_0 \cap CV$ by adding $(n+1)$-cells
was done in a $G$-equivariant way. This implies by passing to the $G$-quotients  the last claim
in  Theorem \ref{thm2}. Alternatively, one can embed $U$ into an affine space $\C^M$ using the Veronese mapping of degree $d_0$ and then use in this new affine setting Proposition 11 in \cite {DP}.
This completes the proof of  Theorem \ref{thm2}.

\endproof

Note also that we have ${\tilde H}^k(U,\Q)={\tilde H}^k(F_0,\Q)^G=0$ for $k<n+1$. In particular ${\tilde H}^k(X,\Q)=0$ for $k<n$, i.e. $X$ is rationally a bouquet of $n$-spheres. In fact $F_0 \cap CV$ can be shown to be a 
bouquet of $n$-spheres and $X=F_0 \cap CV/G$.

\bigskip

{\it Proof of Corollary \ref{cor1}.}

 The first claim follows from Proposition 6.3.6 and Exercise 4.2.13 in  \cite{D} in conjunction to  Theorem 2.10 v  in  \cite{DN}. In fact, to get the vanishing of
$(R^kf_*\C_U)_0$ one has just to write the exact sequence of the triple $(U,T_0,F)$ and to use the fact that
${\tilde H}^k(U,\C)=0$ for $k<n+1$ as we have seen above.
 For the second claim, one has to use Theorem 2.10.i and Proposition 2.18 in  \cite{DN}.
In fact, let $D$ be a large disc in $\C$ containing all the critical values of $f:U \to \C$ inside.
Then $\C^*$ is obtained from $E=\C \setminus D$ by filling in small discs $D_b$ around each critical value
$b \ne 0$ of $f$. In the same way, $M_X$ is obtained from $E_1=f^{-1}(E)$ by filling in the corresponding tubes
$T_b=f^{-1}(D_b)$. It follows from Theorem  \ref{thm4}, (iv), that the inclusion $E_1 \to M_X$ is an 
$n$-equivalence. Now the total space of restriction of the cyclic covering $M_X^c \to M_X$ to the subspace $E_1$
is homotopy equivalent to the generic fiber $F$ of $f$, in such a way that the action of $t$ corresponds to the monodromy at infinity. In this way we get an $n$-equivalence $F \to M_X^c$, inducing the isomorphisms (resp. the monomorphism) announced in Corollary \ref{cor1}, (ii).

To get the similar statement for the monodromy operator $T_0$, we have to build $\C^*$ from a small
punctured disc $D_0^*$ centered at the origin by filling in small discs $D_b$ around each critical value
$b \ne 0$ of $f$. The rest of the above argument applies word for word.

\endproof

The pull-back under $p$ of the infinite cyclic covering $M_X^c \to M_X$ is just
the infinite cyclic covering $( F_0 \setminus CV)   ^c \to   F_0 \setminus CV  $ and we get an induced cyclic covering 
$ p^c:( F_0 \setminus CV)   ^c \to M_X^c  $ 
of order $d_0$. Moreover the action of the deck transformation group $G$
of this covering commutes to the action of the infinite cyclic group $\Z$,
and hence we get the following isomorphism (resp. projection, resp. embedding) of $\Lambda _ {\Q}$-modules

\begin{equation}  \label{eq1}
H^k(M_X^c,\Q)= H_k( ( F_0 \setminus CV)   ^c   ,\Q)_G \leftarrow H_k( ( F_0 \setminus CV)   ^c   ,\Q).
\end{equation}

and

\begin{equation}  \label{eq2}
H^k(M_X^c,\Q)= H^k( ( F_0 \setminus CV)   ^c   ,\Q)^G \to H^k( ( F_0 \setminus CV)   ^c   ,\Q).
\end{equation}

\section{Perverse sheaf approach }

In this section we prove the following weaker version of  Theorem \ref{mainthm}, which is used in the proof
of  Theorem \ref{mainthm}, see subsection 4.2.

\begin{prop}  \label{prop1}

Assume that $d_0$ divides the sum $\sum_{j=1,m}d_j$, say $d d_0=\sum_{j=1,m}d_j$.
Then the  Alexander modules  $H_k(M_X ^c,\C)$ 
  of the hypersurface $X$ are torsion for $k<n+1$. Moreover, 
let $ \lambda \in \C^{*}$ be such that $\lambda  ^d \ne 1$.
Then $ \lambda $ is not a root of the  Alexander polynomials $\Delta _k(t)$ for $m<n+1$.

\end{prop}

The proof we give below to this proposition is closed in spirit to the proofs in \cite{M}, and yields with obvious minor changes (left to the reader)
a proof for our Theorem  \ref{thm3}.

According to  Theorem 4.2 in  \cite{DN}, to prove  Proposition \ref{prop1}, it is enough to prove the following.

\begin{prop}  \label{prop2}

Let $ \lambda \in \C^{*}$ be such that $ \lambda^d \ne 1$, where $d$ is the quotient of $\sum_{j=1,m}d_j$
by $d_0$ .
If $\LL_{\lambda} $ denotes the corresponding local system on $M_W$,
then $H_q(M_W,\LL_{\lambda}  )=0$ for all $q \ne n+1$.

\end{prop}

\proof

First we shall recall the construction of the rank one local system  $\LL_{\lambda} $. Any such local system 
on $M_W$ is given by a homomorphism from $\pi _1(M_W)$ to $\C^*$. To define our local system consider the composition

$$ \pi _1(M_W) \to  \pi _1(M_W') \to  H _1(M_W')= \Z^{m+1}/(d_0,...,d_m) \to \C^*$$
where the first morphism is induced by the inclusion, the second is the passage to the abelianization
and the third one is given by sending the classes $e_0, ...,e_m$ corresponding to the canonical basis of $\Z^{m+1}$ to ${\lambda}^{-d},{\lambda},..., 
{\lambda}$ respectively. For the isomorphism in the middle, see for instance \cite{D0}, p. 102.

\bigskip

It is of course enough to show the vanishing in cohomology, i.e.
 $H^q(M_W,\LL_{\lambda}  )=0$   for all $q \ne n+1$. Let $i:M_W \to  U$ and $j: U \to Z$ be the two inclusions. Then one clearly has
$ \LL_{\lambda}   [n+1] \in Perv(M_W)$ and hence $\F=Ri_*(  \LL_{\lambda}  [n+1]) \in Perv( U   )$,
since the inclusion $i$ is a quasi-finite affine morphism. See for this and the following p. 214 in  \cite{D} for a similar argument.

Our vanishing result will follow from a study of the natural morphism
$$ Rj_!\F \to  Rj_*\F.$$ 
Extend it to a distinguished triangle
$$ Rj_!\F \to  Rj_*\F   \to \G  \to  . $$
Using the long exact sequence of hypercohomology coming from the above triangle, we see exactly as on p.214 in \cite{D} that all we have to show is that
$\HH^k(Z,\G)=0$ for all $k<0$. This vanishing obviously holds if we show that $\G=0$.

This in turn is equivalent to the vanishing of all the local cohomology groups of $ Rj_*\F$, namely $H^m(M_x,\LL_x)=0$ for all $m \in \Z$ and for all points $x \in W_0$.
Here  $M_x = M_W \cap B_x$ for $B_x$ a small open ball at $x$ in $ Z$ and $\LL_x$ is the restriction of the local system $\LL_{\lambda}  $ to $M_x$.

The key observation is that, as already stated above, the action of an oriented elementary loop about the hypersurface $W_0$ in the  local systems  $\LL_{\lambda}  $ and $\LL_x$ corresponds to multiplication by 
$\nu = \lambda ^{-d} \ne 1$.

There are two cases to consider.

\bigskip

\noindent{Case 1.} If $x \in W_0 \setminus V $, then $M_x$ is homotopy equivalent to
 $\C^*$ and the corresponding local system $\LL_{\nu}  $ on  $\C^*$ is defined by multiplication by $\nu$, hence the claimed vanishings are obvious.

\bigskip

\noindent{Case 2.} If  $x \in W_0 \cap V $, then due to the local product structure of stratified sets cut by a transversal, $M_x$ is homotopy equivalent to a product $(B' \setminus (V \cap B')) \times \C^*$, with $B'$ a small open ball centered at $x$ in $W_0$, and the corresponding local system is an external tensor product, the second factor being exactly $\LL_{\nu}    $.  The claimed vanishings follow then from the K\"unneth Theorem, see 4.3.14 \cite{D}.

\endproof

A minor variation of this proof gives also  Theorem \ref{thm3}. Indeed, let $D=\sum _{j=0,m}d_j$ and let 
$\alpha $ be a $D$-root of unity, $\alpha \ne 1$. All we have to show is that  $H^q(M_W,\LL_{\alpha}  )=0$   for all $q \ne n+1$, see for instance 6.4.6 in  \cite{D}.

The action of an oriented elementary loop about the hypersurface $W_0$ in the  local systems  $\LL_{\alpha}  $ and in its restrictions $\LL_x$ as above corresponds to multiplication by 
$\alpha \ne 1$. Therefore the above proof works word for word.

\bigskip

One has also the following result, in which the bounds are weaker than in Maxim's Theorem 4.2 in \cite{M}.

\begin{prop}  \label{prop3}

Assume that $d_0$ divides the sum $\sum_{j=1,m}d_j$, say $d d_0=\sum_{j=1,m}d_j$.
Let $ \lambda \in \C^{*}$ be such that $\lambda ^d = 1$ and let $\sigma$ be a non negative integer.
Assume that $ \lambda   $ is not a root of the $q$-th  local Alexander polynomial $\Delta _q (t)_x$ of the hypersurface singularity $(V,x)$ for any $q <n+1 - \sigma$ and any point $x \in W_1$, where $W_1$ is an irreducible component of $W$ different from $W_0$.

Then $ \lambda   $ is not a root of the global Alexander polynomials   $\Delta _q (t)$ associated to $X$ for any
 $q <n+1 - \sigma$.

\end{prop}

\bigskip

To prove this result, we start by the following general remark.

\begin{remark}\label{rm1}

 If $S$ is an $s$-dimensional stratum in a Whitney stratification of $V$ such that $x \in S$ and $W_0$ is transversal to $V$ at $x$, then, due to the local product structure,  the $q$-th reduced local Alexander polynomial $\Delta _q (t)_x$ is the same as that of the  hypersurface singularity $V \cap T$ obtained by cutting
the germ $(V,x)$ by an $(n+1-s)$-dimensional transversal $T$. It follows that these reduced local Alexander polynomials $\Delta _q (t)_x$ are all trivial except for $q \leq n-s$. It is a standard fact that, in the local situation of a hypersurface singularity, the  Alexander polynomials can be defined either from the link or as the characteristic polynomials of the corresponding  the monodromy operators. Indeed, the local Milnor fiber is homotopy equivalent to the corresponding infinite cyclic covering.

\end{remark}

 Let $i:M_W \to  Z \setminus W_1  $ and $j:  Z \setminus W_1      \to Z$ be the two inclusions. Then one  has
$ \LL_{\lambda}    [n+1] \in Perv(M_W)$ and hence $\F=Ri_*( \LL_{\lambda}   [n+1]) \in Perv( Z \setminus W_1    )$,
exactly as above.

Extend now the natural morphism $ Rj_!\F \to  Rj_*\F$ to a distinguished triangle
$$ Rj_!\F \to  Rj_*\F   \to \G .   $$

Applying Theorem 6.4.13 in  \cite{D} to this situation, and recalling the above use of  Theorem 4.2 in \cite{DN}, all we have to check is that $H^m(M_x,\LL_x)=0$ for all points $x \in W_1$  and $m<n+1 - \sigma$.  For $ x \in W_1 \setminus W_0$, this claim is clear by the assumptions made. The case when $x \in W_1 \cap W_0$  can be treated exactly as above, using the product structure, and the fact that
the monodromy of $(W_1,x)$ is essentially the same as that of  $(W_1 \cap W_0 ,x)$, see our remark above.

This completes the proof of  Proposition \ref{prop3}.

\begin{remark}\label{rm2}
Here is an alternative explanation for some of the bounds given in  Theorem 4.2 in \cite{M}.
Assume that $ \lambda $ is a root of the Alexander polynomial 
 $\Delta _i (t)$ for some $i<n+1$. Then it follows from Proposition \ref{prop3}
the existence of a point $x \in W_1$ and of an integer $\ell \leq i$ such that
 $\lambda  $ is a root of the local Alexander polynomial $\Delta _{\ell} (t)_x$.
If $x \in S$, with $S$ a stratum of dimension $s$, then by Remark \ref{rm1}, we have $\ell \leq n-s$. This provides half of the bounds in  Theorem 4.2 in \cite{M}. The other half comes from the following remark. Since 
 $\lambda  $ is a root of the Alexander polynomial 
 $\Delta _i (t)$, it follows that $H^i(M_W, \LL_{\lambda}      ) \ne 0$. This implies via an obvious exact sequence that $\HH^{i-n-1}(W_1,\G)\ne 0$. Using the standard spectral sequence to compute this hypercohomology group, we get that some of the groups
$H^p(W_1,\HC^{i-n-1-p}\G)$ are non zero. This can hold only if $p \leq 2 dim (Supp\HC^{i-n-1-p}\G)$. Since $\HC^{i-n-1-p}\G_x = H^{i-p}(M_x,L_x)$ his yields the inequality $p=i- \ell \leq 2s$ in
 Theorem 4.2 in \cite{M}. 
\end{remark}

\begin{remark}\label{rm3}

Let $ \lambda \in \C^{*}$ be such that $ \lambda^d = 1$, where $d$, the quotient of $\sum_{j=1,m}d_j$
by $d_0$, is assumed to be an integer.
Let $\LL_{\lambda} $ denotes the corresponding local system on $M_W$. The fact that the associated monodromy about the divisor $W_0$ is trivial can be restated as follows. Let  $\LL_{\lambda}' $ be the rank one local system on $M_V=Z \setminus V$ associated to ${\lambda}$. Let $j:M_W \to M_V$ be the inclusion. Then
$$\LL_{\lambda}=j^{-1}\LL_{\lambda}'.$$
Let moreover  $\LL_{\lambda}'' $ denote the restriction to  $\LL_{\lambda}' $ to the smooth divisor 
$W_0 \setminus (V \cap W_0)$. Then we have the following Gysin-type long exact sequence
$$ ...  \to  H^q(M_V,\LL_{\lambda}') \to  H^q(M_W,\LL_{\lambda}) 
\to H^{q-1}(W_0 \setminus (V \cap W_0),\LL_{\lambda}'')  \to  H^{q+1}(M_V,\LL_{\lambda}')    \to    ...$$
exactly as in  \cite{D} , p.222.

 The cohomology groups $ H^*(M_V,\LL_{\lambda}')$ and $ H^*(W_0 \setminus (V \cap W_0),\LL_{\lambda}'')$ being usually simpler to compute than $ H^*(M_W,\LL_{\lambda})$, this exact sequence can give valuable information on the latter cohomology groups.

\end{remark}

\section{Semisiplicity results}

In this section we prove the first claim in our main result Theorem \ref{mainthm}.

\subsection{ First proof (the case  $W_0=H$ is the hyperplane at infinity in $\PP^{n+1}$.)}

\ \\
 
Let $\U$ be a sufficiently small 
tubular neighborhood of the hyperplane $H$ at infinity.
We claim the following: 
\bigskip 

\noindent (i) 
$ \ \ 
 \pi_i(\U \setminus (H \cup (V\cap \U))) \rightarrow \pi_i(M_X) \ \ {\rm  is \ an \ isomorphism \ for}
\ 1 \le i \le n-1$

\bigskip
\noindent (ii)
 $ \ \ \pi_n(\U \setminus (H \cup (V\cap \U))) \rightarrow \pi_n(M_X) \ \ {\rm is \ \ surjective}.$

\bigskip

First notice that as a consequence of transversality of $V$ and $H$ 
we have $S^1$-fibration: 
$\U \setminus (H \cup (V \cap \U)) \rightarrow H \setminus (H \cap V)$.
Indeed, if $f(x_0,...,x_{n+1})=0$ is an equation of 
$V$ and $x_0=0$ is the equation for $H$ then 
the pencil $\lambda f(x_0,...,x_{n+1})+ \mu f(0,x_1,...,x_{n+1})$
defines deformation of $V$ to the cone over $V \cap H$. 
Since $V$ is transversal to $H$ this pencil contains   
isotopy of $\U \cap V$ into the intersection of $\U$ with 
the cone.

Let $Y$ denotes the above cone in 
$\PP^{n+1}$ over $V \cap H$. The obvious  
$\C^*$-bundle $\PP^{n+1} \setminus (Y \cup H)
 \rightarrow H \setminus (H \cap V)$ is homotopy equivalent to 
the above $S^1$-bundle: 
$\U \setminus (H \cup (V \cap \U)) \rightarrow H \setminus (H \cap V)$.
We can apply to both $M_X$ and $\PP^{n+1} \setminus (Y \cup H)$ the 
Lefschetz hyperplane section theorem for stratified spaces
(cf. \cite{GM}, theorem 4.3) using a generic hyperplane $H'$. 
Thus for  $i \le n-1$ we obtain the isomorphisms: 
$$\pi_i(M_X)=\pi_i(M_X \cap H')=
\pi_i((\PP^{n+1} \setminus (Y \cup H)) \cap H')=\pi_i(\PP^{n+1} 
\setminus (Y \cup H))$$
(the middle isomorphism takes place since for $H'$ near $H$ 
both spaces are isotopic). This yields (i).

To see (ii), let us apply Lefschetz hyperplane section theorem
to a hyperplane $H'$ belonging to $\U$. We obtain the surjectivity 
of the map which is the following composition:
$$
\pi_i(H' \setminus (V \cup H)) \rightarrow \pi_i(\U \setminus (H' \cup 
(V \cap \U))) \rightarrow \pi_i(M_X) 
$$
Hence the right map is surjective as well.

\bigskip
The relations (i) and (ii)
yield that  $M_X$ has the homotopy type of 
a complex obtained 
from $\U \setminus (H \cup (V\cap \U))$ by adding cells having 
the dimension greater than or equal to $n+1$.
 Hence the same is true for the infinite 
cyclic covers defined as in section 1 for
 $M_X$ and $\U \setminus (H \cup (V \cap \U))$ 
respectively. Denoting by $(\U \setminus (H \cup (V \cap \U)))^c$
the infinite cyclic cover of the latter we obtain that

\begin{equation}\label{surjectivity}
H_i({(\U \setminus \U \cap (V \cup H))}^c,\Q) 
\rightarrow H_i(M_X^c,\Q )
\end{equation}
 is surjection for $i=n$ and the isomorphism for $i < n$.
Since the maps above are induced by an embedding map, 
they are isomorphisms or surjections of $ \Lambda_{\Q}  $-modules.

\par As was mentioned above, since $V$ is transversal to $H$, the space 
$\U \setminus (\U \cap (V \cup H))$ 
is homotopy equivalent to the complement in affine space 
to the cone over the projective hypersurface $V \cap H$. 
On the other hand, the complement in $\C^{n+1}$ to the cone over $V \cap H$  
is homotopy equivalent to the 
complement to $V \cap S^{2n+1}$ in $S^{2n+1}$ 
where $S^{2n+1}$ is a sphere about the vertex of the cone.
The latter, by the Milnor's theorem (cf. \cite{milnor1})
is fibered over the circle. 
Hence the fiber of this fibration, as the Milnor fiber 
of any hypersurface singularity, is homotopy equivalent 
to the infinite cyclic cover of $S^{2n+1} \setminus V \cap S^{2n+1} 
\approx \C^{n+1} \setminus V$.
As in section 1, this cyclic cover is the one 
corresponding to the kernel of the homomorphism 
of the fundamental group given by the linking number.
In particular, since a Milnor fiber is a finite CW-complex, 
$H_i({\U   \setminus   (\U \cap (V \cup H))}^c,\Q)$
is a finitely generated $\Q$-module and hence a torsion
$\Lambda_{\Q}$-module.
Moreover, the homology of 
the Milnor fiber of a cone 
and hence $H_i({\U  \setminus  (\U \cap (V \cup H))}^c,{\C})$
 is annihilated by $t^d-1$ since 
the monodromy on $f(0,x_1,...,x_{n+1})=1$ 
is given by multiplication of coordinates by a root of 
unity of degree $d$ and hence has the order equal to $d$.
Therefore it follows from the surjectivity of (\ref{surjectivity})
 that  the same is true for $H_i(M_X^c)$.
In particular $H_i(M_X^c)$ is semisimple.

\bigskip 

\subsection {Second proof  (the general case)}

\   \\

Using the equation \ref{eq1} and Proposition \ref{prop1}, it is enough to show that the Alexander invariant 
$A_k= H_k( ( F_0 \setminus CV)   ^c   ,\Q)$ of the hypersurface $h=0$ in the affine variety $F_0$
is a torsion semisimple $\Lambda _{\Q}$-module killed by $t^e-1$ for some integer $e$.
Indeed, one we know that $t$ is semisimple on $H_k(M_X^c, \Q)$, Proposition \ref{prop1} implies that $t^d=1$.

The fact that $A_k$ is torsion follows from Theorem 2.10.v in \cite{DN} and the claim (iv) in Theorem  \ref{thm4}.
Moreover, Theorem 2.10.ii in \cite{DN} gives for $k \le n$, an epimorphism of $\Lambda _{\Q}$-modules
\begin{equation}\label{eq3}
H_k(F_1,\Q) \to A_k
\end{equation}
where $F_1$ is the generic fiber of $h:F_0 \to\C$ and $t$ acts on $H_k(F_1,\Q)$ via the monodromy at infinity.
By definition, the monodromy at infinity of $h:F_0 \to\C$ is the monodromy of the fibration over the circle $S^1_R$,
 centered at the origin and  of radius
$R>>0$, given by 
$$\{x \in CZ;~~f_0(x)=1,~~|h(x)|=R\} \to  S^1_R, ~~~~ x \mapsto h(x).$$
Using a rescaling, this is the same as the fibration
\begin{equation}\label{eq4}
\{x \in CZ;~~f_0(x)= \epsilon ,~~|h(x)|= 1\} \to  S^1, ~~~~ x \mapsto h(x)
\end{equation}
where $0< \epsilon <<1$.

Let $R_1>>0$ be such that
$$\{x \in CZ;~ |x| \le R_1 ,~~|h(x)|= 1\} \to  S^1, ~~~~ x \mapsto h(x)$$
is a proper model of the Milnor fibration of $h:CZ \to\C$. This implies that all the fibers
$\{h=t\}$ for $t\in S^1$ are transversal to the link $K=CZ \cap S_{R_1}^{2N+1}.$

A similar argument, involving the Milnor fibration of  $h:CW_0 \to \C$ shows that
 all the fibers
$\{h=t\}$ for $t\in S^1$ are transversal to the link $K_0=CW_0 \cap S_{R_1}^{2N+1}$.
Using the usual $S^1$-actions on these two links, we see that transversality for all fibers
$\{h=t\}$ for $t\in S^1$ is the same as transversality for $\{h=1\}$. But saying that $\{h=1\}$
is transversal to $K_0$ is the same as saying that $Z_{0,1} \pitchfork K$. By the compactness of $K$,
there is a $\delta >0$ such that  $Z_{s,1} \pitchfork K$ for $|s|<\delta$. Using the above
 $S^1$-actions on  links, this implies that  $Z_{s,1} \pitchfork K$ for $|s|<\delta$ and  $t\in S^1$.

Choose $\delta $ small enough such that the open disc $D_{\delta}$ centered at the origin and  of radius
 $\delta $ is disjoint from the finite set of circles $g_0(\Gamma (h,g_0) \cap h^{-1}(S^1))$.
Using the relative Ehresmann Fibration Theorem, see for instance  \cite{D0}, p.15, we see that the map
$$\{x \in CZ;~|x| \le R_1, ~~f_0(x)< \delta ,~~|h(x)|= 1\} \to D_{\delta}\times S^1   , ~~~~ x \mapsto (g_0(x),h(x))$$
is a locally trivial fibration. It follows that the two fibrations
$$\{x \in CZ;~|x| \le R_1, ~~f_0(x)= \delta/2 ,~~|h(x)|= 1\} \to S^1   , ~~~~ x \mapsto h(x)$$
and
$$\{x \in CZ;~|x| \le R_1, ~~f_0(x)=0 ,~~|h(x)|= 1\} \to S^1   , ~~~~ x \mapsto h(x)$$
are fiber equivalent. In particular they have the same monodromy operators.
The first of these two fibration is clearly equivalent to the  monodromy at infinity fibration  \ref{eq4}.
The homogeneity of the second of these two fibration implies that its monodromy operator has order
$e=dd_0$. This ends the proof of the semisimplicity claim in the general case.

\section{ Mixed Hodge structures on Alexander invariants}

This proof involves several mappings and the reader may find useful to draw them all in a diagram.

Since the mapping $f:M_X \to \C^*$ has a monodromy of order $d$ (at least in dimensions $k<n$, see Corollary
\ref{cor1} and Theorem \ref{mainthm}, (i)), it is natural to consider the base change $\phi :  \C^* \to  \C^*$ given by $s \mapsto s^d$. Let $f_1: M_X^d \to \C^*$ be the pull-back of $f:M_X \to \C^*$ under $\phi$ and
let  $\phi _1:  M_X^d \to M_X$ be the induced mapping, which is clearly a cyclic $d$-fold covering.
It follows that the infinite cyclic covering $p_c: M_X^c \to M_X$ factors through $ M_X^d$, i.e. there is
an infinite cyclic covering $p_d: M_X^c \to M_X^d$ corresponding to the subgroup $<t^d>$ in  $<t>$, such that
$\phi _1 \circ p_d =p_c$. Since $ M_X^d= M_X^c/<t^d>$, it follows that $t$ induces an automorphism $t$ of
$ M_X^d$ of order $d$.

Let $F_1= f_1^{-1}(s)$ be a generic fiber of $f_1$, for $|s|>>0$. Then $\phi _1$ induces a
regular homeomorphism $F_1 \to F=f^{-1}(s^d)$. Let $i:F \to M_X$ and $i_1:F_1 \to M_X^d$ be the two inclusions.
Note the $i$ has a lifting $i_c:F \to  M_X^c$, which is exactly the $n$-equivalence mentionned in the proof
of the Corollary \ref{cor1}, commuting at the cohomology level with the actions of $t$ and $T_{\infty}$.
Moreover, $i$ has a lifting $i_d=p_d \circ i_c$ such that $i=\phi _1 \circ i_d$.

Now we consider the induced morphisms on the various cohomology groups.
It follows from the general spectral sequences relating the cohomology of  $M_X^c$ and $M_X^c/<t^d>$,
see \cite{W}, p. 206, that $p_d^*:H^k(M_X^d) \to H^k(M_X^c)$ is surjective. It follows that
$H^k(M_X^c)$ is isomorphic (as a $\Q$-vector space endowed to the automorphism $t$) via $i_c^*$ to the
sub MHS in $H^k(F)$ given by $i_d^*(H^k(M_X^d))$. Note that $i_d$ can be realized by a regular mapping
and $i_d^*$ commutes with the actions of $t$ and  $T_{\infty}$.

There is still one problem to solve, namely to show that this MHS is independent of $s$, unlike the MHS
 $H^k(F,\Q)$ which depends in general on $s$, see the example below.
To do this, note that $\phi _1^* \circ i_d^*(H^k(M_X^d))=i_1^*(H^k(M_X^d))$ as MH substructures in 
$H^k(F_1,\Q)$. More precisely, $i_1^*(H^k(M_X^d))$ is contained in the subspace of invariant cocycles 
$H^k(F_1,\Q)^{inv}$, where ${inv}$ means invariant with respect to the monodromy of the mapping $f_2:M_2 \to S_2$
oblained from $f_1$ by deleting all the singular fibers, e.g. $S_2=\C^* \setminus C(f_1)$, where $C(f_1)$
is the finite set of critical values of $f_1$. We have natural morphisms of MHS
$$H^k(M_X^d) \to H^k(M_2) \to H^0(S_2,R^kf_{2,*}\Q)$$
the first induced by the obvious inclusion, the second coming from the Leray spectral sequence
of the map $f_2$, see  \cite{Sa1} (5.2.17-18), \cite{Sa2} (4.6.2) and \cite{Sa3}. Moreover the last morphism above is surjective.
On the other hand, there is an isomorphisms of MHS
$$H^0(S_2,R^kf_{2,*}\Q) \to H^k(F_1,\Q)^{inv}$$
showing that the latter MHS is independent of $s$, see  \cite{SZ}, Prop. (4.19). It follows that
$i_1^*(H^k(M_X^d))$ has a MHS which is independent of $s$. By transport, we get a natural MHS on 
$H^k(M_X^c,\Q)$ which clearly satisfies all the claim in Theorem \ref{mainthm}, (ii).

\begin{example} \label{ex2} \rm

For $f:\C^2 \to \C$ given by $f(x,y)=x^3+y^3+xy$, let $F_s$ denote the fiber $f^{-1}(s)$. Then, the MHS on
$H^1(F_s,\Q)$ (for  $F_s$ smooth) depends on $s$. Indeed, it is easy to see that the graded piece
$Gr^W_1H^1(F_s,\Q)$ coincides as a Hodge structure to $H^1(C_s,\Q)$, where $C_s$ is the elliptic curve
$$x^3+y^3+xyz-sz^3=0.$$
Moreover, it is known that $H^1(C_s,\Q)$ and $H^1(C_t,\Q)$ are isomorphic as Hodge structures if and only if
the elliptic curves $C_s$ and  $C_t$ are isomorphic, i.e. $j(s)=j(t)$, where $j$ is the $j$-invariant of an elliptic curve. This proves our claim and shows that the range in Corollary  \ref{cor2} is optimal.

For $f:\C^2 \to \C$ given by $f(x,y)=(x+y)^3+x^2y^2$ it is known that the monodromy at infinity operator has a Jordan block of size 2 corresponding to the eigenvalue $\lambda=-1$, see  \cite{GN}.

\end{example}

\begin{remark} A ``down to earth'' relation between the cohomology 
of $M_X^c$ and $M_X^d$ used above and obtained from the spectral sequence
\cite{W} can be described also using the ``Milnor's exact sequence''
i.e. the cohomology sequence corresponding to the 
sequence of chain complexes:
$$0 \rightarrow C_*(M_X^c) \rightarrow C_*(M_X^c) \rightarrow C_*(M_X^d)
\rightarrow 0$$
This is a sequence of free $\Q[t,t^{-1}]$-modules with the 
left homomorphism given by multiplication by $t^d-1$.
The corresponding cohomology sequence is:

\begin{equation}\label{milnorseq}
0 \rightarrow 
H^{i}(M_X^c) \buildrel  \iota \over
\rightarrow H^{i+1} (M_X^d) \rightarrow H^{i+1}(M_X^c) 
\rightarrow 0
\end{equation}

The zeros on the left and the right in (\ref{milnorseq}) appear
because of mentioned earlier triviality of the action of $t^d$
on cohomology. 
Another way to derive (\ref{milnorseq}) is to consider the 
the Leray spectral sequence corresponding to the classifying map 
$M^c_X \rightarrow BS^1=\C^*$ corresponding to the action of $t$.
This spectral sequence degenerates in term $E_2$ and is equivalent to 
the sequence (\ref{milnorseq}).
A direct argument shows that the image of $\iota$ coincides with the 
kernel of the cup product $H^1(M_X^d) \otimes H^{i+1}(M_X^d) \rightarrow 
H^{i+2}(M_X^d)$ (i.e. the annihilator of $H^1(M_X^d)$)
which also yields the MHS on 
$H^{i}(M_X^c)$ as a subMHS on $H^{i+1}(M_X^d)$.
\end{remark}

\begin{remark}\label{superab} 

The above Mixed Hodge structure plays  a key role 
in the calculation of the first non-vanishing 
homotopy group of the complements 
to a hypersurface $V$  in $\PP^{n+1}$ 
with isolated singularities (cf. \cite{Lib1}). 
More precisely, in this paper for each 
$\kappa=exp {{2 \pi k \sqrt {-1}}\over d}$
(or equivalently $k, 0 \le k \le d-1$)
and each point $P \in V \subset \PP^{n+1}$ 
which is singular on $V$ 
the ideal 
${\mathcal A}_{ P , \kappa } $ 
is associated
(called there {\it the ideal of quasiadjuncton}).
These ideals glued together into a subsheaf 
${\mathcal A}_{\kappa}  \subset {\mathcal O}_{\PP^{n+1}}$ 
of ideals having at $P$ the stalk 
${\mathcal A}_{P,\kappa }$ and ${\mathcal O}_Q$ and any 
other $Q \in \PP^{n+1} \setminus Sing(V)$. It is shown in \cite{Lib1} that 
for the $\kappa$-eigenspace
of $t$ acting on $F^0H^n(({\PP^{n+1}  \setminus     (V \cup H)})^c)$ one has:
$$F^0H^n(( {\PP^{n+1}  \setminus     (V \cup H)}        )^c)_{\kappa}=
{\rm dim}H^1({\mathcal A}_{\kappa} (d-n-2-k))$$
The right hand side can be viewed as the difference between 
actual and ``expected'' dimensions of the linear system of hypersurfaces 
of degree $d-n-2-k$ which local equations belong 
to the ideals of quasiadjunction at the singular points of $V$.  
In the case of plane curves see also  \cite{Es}, \cite{LV}.

\end{remark}

\end{document}